\documentclass[10pt,conference]{IEEEtran}
\usepackage{amssymb}
\usepackage{amsmath}
\usepackage{amsfonts}
\usepackage{verbatim}
\usepackage{color}
\usepackage{multicol}

\newtheorem{corollary}{Corollary}

\newtheorem{proposition}{Proposition}

\newtheorem{theorem}{Theorem}

\begin{document}
\title{A planar large sieve and sparsity of  time-frequency representations }
\author{\IEEEauthorblockN{Lu\'{\i}s Daniel Abreu}
\IEEEauthorblockA{ Acoustics Research Institute,\\
Wohllebengasse 12-14, Vienna A-1040, Austria.\\
Email: labreu@kfs.oeaw.ac.at}
\and
\IEEEauthorblockN{Michael Speckbacher}
\IEEEauthorblockA{ Acoustics Research Institute,\\
Wohllebengasse 12-14, Vienna A-1040, Austria.\\
Email: speckbacher@kfs.oeaw.ac.at}
}
\maketitle

\textbf{\textit{Abstract -- }
With the aim of measuring the sparsity of a real signal, Donoho and Logan
introduced the concept of maximum Nyquist density, and used it to extend
Bombieri's principle of the large sieve to bandlimited functions. This led
to several recovery algorithms based on the minimization of the $\mathbf{L_{1}}$%
-norm. In this paper we introduce the concept of {\ planar maximum} Nyquist
density, which measures the sparsity of  the time-frequency distribution
of a function. We obtain a planar large sieve principle which applies  to time-frequency representations with a gaussian window, or equivalently,
to Fock spaces, 
allowing for %
perfect recovery of the short-Fourier transform (STFT) of functions in the modulation space $\mathbf{M_{1}}$ (also known as Feichtinger's algebra $\mathbf{S_{0}}$) corrupted by sparse
noise and for approximation of missing STFT data in $\mathbf{M_{1}}$, by
$\mathbf{L_{1}}$-minimization.
}

\section{Introduction}

With the aim of measuring the sparsity of a real signal, Donoho and Logan
introduced the concept of \emph{maximum Nyquist density}, defined in \cite%
{DonohoLogan} as 
\begin{align*}
\rho (T,W)&:=W\cdot \sup_{t\in \mathbb{R}}\left\vert T\cap \lbrack
t,t+1/W]\right\vert\\ &\leq W\cdot \left\vert T\right\vert \text{,}
\end{align*}%
where $T\subset \mathbb{R}$ and $W$ is the band-size in the space of
band-limited functions 
\begin{equation*}
B_W:=\left\{ f\in L^{1}:\ \hat{f}(\omega)=0,\ \forall\ |\omega|> \pi W\right\} \text{.}
\end{equation*}%
If the set $T$ is sparse in terms of low Lebesgue measure (small concentration in any
interval of length $1/W$), then $\rho (T,W)$ can be
considerably small compared to the natural Nyquist density $W\cdot
\left\vert T\right\vert $. We will write $P_{A}f=\chi _{A}f$ for the
multiplication by the indicator function of $A$. In \cite[Theorem 7]%
{DonohoLogan}, Donoho and Logan proved that, if $f\in B_W$ and $\delta <\frac{2 }{W}$,
the inequality%
\begin{equation}\label{DLl1}
\Vert P_{T}f\Vert _{1}
\leq \sup_{t\in 
\mathbb{R}}\left\vert T\cap \lbrack t,t+\theta ]\right\vert \cdot\frac{\pi W/2}{\sin \left(\pi W\theta /2\right)}
\cdot\Vert f\Vert_1
\end{equation}
holds. In particular, if $\delta(T)$ denotes the norm of the projection operator $P_T$ and $\theta=1/W$, then $\delta(T)\leq\frac{\pi}{2}\rho(T,W)$. Note that the inequality (\ref{DLl1}) 
falls within the realm of quantitative uncertainty principles \cite{DonohoStark,JamPow,JamPowJFA,Eugenia,Tao}, which paved the
way to the modern theory of compressed sensing (see \cite[Section 1.6]{CRT} or \cite%
{Donoho,FR}).

Donoho and Logan's
interest in such inequalities, in particular in obtaining good constants
depending on the sparsity of the set $T$, was motivated by signal recovery
problems \cite{DonohoStark}. As an application, they derived the
following results, which allows to perfectly reconstruct a bandlimited signal
corrupted by sparse noise using $L_{1}$-norm minimization.

\textit{Corollary:} (\cite[Corolllary 1]{DonohoLogan}) Suppose that $g=b+n$ {is
observed}, $b\in B_W$, $n\in L_1$, and that the{\ }\emph{unknown support} $T$ of the
noise $n$ satisfies 
\begin{equation}
\rho (T,W)<1/\pi \text{.}  \label{sparse}
\end{equation}%
Then the solution of the minimization problem 
\begin{equation*}
\beta (g)=\arg \min_{\widetilde{b}\in B_W}\big\Vert g-\widetilde{b}%
\big\Vert_{1}
\end{equation*}%
is unique and recovers the signal $b$ perfectly ($\beta(g)=b$).
\medskip

This extends the so-called Logan's phenomenon (\cite{Logan}, see also the
discussion in \cite[Section 6.2]{DonohoStark}). 
The following recovery result for missing data extends the results from \cite{DonohoStark} from the $L_2$ to the $L_1$ setting. We will give a proof in the STFT context in Corollary \ref{missingcorollary}.

\textit{Corollary:} Suppose that we observe $h=P_{T^c}(b+n)$, where $b\in B_W$, $\|n\|\leq\varepsilon$ and 
\begin{equation*}
\rho(T,W)<2/\pi. 
\end{equation*}
Then, any solution of the minimization problem 
$$\sigma(h):=\text{arg}\min_{\widetilde b\in B_W}\|P_{T^c}(\widetilde b-h)\|_1$$ satisfies 
$$
\|b-\sigma(h)\|_1\leq \frac{4\varepsilon}{2-\pi\cdot \rho(T,W)}.
$$
\smallskip

 Let $\varphi (t)=2^{1/4}e^{-\pi t^{2}}$ be the normalized
gaussian. The \emph{short-time Fourier transform (STFT)} is defined as
follows 
\begin{equation}
V_{\varphi }f(x,\omega )=\int_{%
\mathbb{R}
}f(t)\overline{\varphi (t-x)}e^{-2\pi i\omega t}dt\text{.}  \label{Gabor}
\end{equation}%
Moreover, define the modulation spaces 
\begin{equation*}
M^{p}:=\big\{f\in \mathcal{S}^{\prime }(\mathbb{R}):\ \Vert V_{\varphi
}f\Vert _{p}<\infty \big\},\ \ p\geq 1.
\end{equation*}%

Modulation spaces are ubiquitous in time-frequency analysis \cite{FG,Charly}%
. They were introduced in \cite{FeiModulation}. It is the
purpose of this paper to obtain a planar version of (\ref{DLl1}) and apply
it to recovery  problems for the short-time Fourier transform of functions in $M^{1}$ using $L_{1}$-minimization. The space $M^1$
 is also known as Feichtinger's algebra $S_{0}$ and it can be
identified with the Bargmann-Fock space $\mathcal{F}_{1}(\mathbb{C})$ of
entire functions.

As a planar analogue of $\rho (T,W)$, we introduce the following concept.
Let $\Delta \subset \mathbb{R}^{2}$. The planar maximum Nyquist density $%
\rho (\Delta ,R)$ is defined as 
\begin{equation}
\rho (\Delta ,R):=\sup_{z\in \mathbb{R}^{2}}|\Delta \cap (z+D_{1/R})|\leq
|\Delta |\text{,}  \label{planarNyquist}
\end{equation}%
where $D_{1/R}\subset \mathbb{R}^{2}$ is the disc of radius $1/R$ centered
in the origin. If the set $\Delta $ is sparse in the sense of Lesbegue
measure (small concentration in any disc of radius $1/R$) then $\rho (\Delta ,R)$\ can be
considerably smaller than the natural Nyquist density\ $|\Delta |$ (see %
\cite{AGR,Daubechies,DeMarie,Ly,RS95,Seip0} for natural Nyquist densities
in the context of Fock and modulation spaces and \cite{FGHKR} for a survey
on the current state of the art of the topic). Our main result is the
following.

\begin{theorem}\label{main}
Consider  $\Delta \subset \mathbb{R}^{2}$ and let $f\in M^{1}$, then,
for every $0<R<\infty $, it holds
\begin{equation}
\Vert P_{\Delta }(V_{\varphi }f)\Vert _{1}\leq \frac{\rho (\Delta ,R)}{%
1-e^{-\pi /R^{2}}}\cdot \Vert V_{\varphi }f\Vert _{1}\text{.}  \
\end{equation}
\end{theorem}
\medskip 

Set 
\begin{equation*}
\delta (\Delta ):=\sup_{f\in M^{1}}\dfrac{\Vert P_{\Delta }(V_{\varphi
}f)\Vert _{1}}{\Vert V_{\varphi }f\Vert _{1}}\text{.}
\end{equation*}%
By Theorem 1, 
\begin{equation}\label{delta-bound}
\delta (\Delta )\leq (1-e^{-\pi /R^{2}})^{-1}\rho (\Delta ,R).
\end{equation}
 Moreover, if $\delta (\Delta )<\frac{1}{2}$ then every $f\in M^{1}$
satisfies $\Vert P_{\Delta }(V_{\varphi }f)\Vert _{1}<\Vert P_{\Delta
^{c}}(V_{\varphi }f)\Vert _{1}$. Combined with the argument in \cite[Section
6.2]{DonohoStark}, this implies the following result, which allows to
perfectly reconstruct the STFT of a signal in $M^1$
corrupted by sparse noise, using $L_{1}$-norm minimization.

\begin{corollary}
Suppose that $G=V_{\varphi }f+N$ is observed, where $f\in M^{1}$, $N\in L_{1}(%
\mathbb{R}^{2})$ and that the \emph{unknown support} $\Delta $ of $N$
satisfies 
\begin{equation}
\rho (\Delta ,R)<\frac{1}{2}(1-e^{-\pi /R^{2}})\text{,}  \label{planarsparse}
\end{equation}%
for some $R>0$.
Then $\delta (\Delta )<\frac{1}{2}$ and the solution of the minimization
problem 
\begin{equation*}
\beta (G)=\arg \min_{g\in M^{1}}\big\Vert G-V_{\varphi }g\big\Vert_{1}
\end{equation*}%
is unique and recovers the signal $f$ perfectly ($\beta
(G)=f$).
\end{corollary}
\medskip

One can also derive an analogue for the recovery of missing data.
\begin{corollary}\label{missingcorollary}
Let $f\in M^1$ and suppose that one observes $H=P_{\Delta^c}(V_\varphi f+N)$, where $\|N\|_1\leq \varepsilon$ and that the domain $\Delta$ of missing data satisfies
\begin{equation}\label{planarsparselessone}
\rho (\Delta ,R)<(1-e^{-\pi /R^{2}})\text{,} 
\end{equation}%
for some $R>0$. Then any solution of 
\begin{equation*}
\sigma (H)=\arg \min_{h\in M^{1}}\big\Vert P_{\Delta^c}( H-V_{\varphi }h)\big\Vert_{1}
\end{equation*}
satisfies
$$
\left\Vert V_\varphi \big(f-\sigma(H)\big)\right\Vert_1\leq\frac{2\varepsilon(1-e^{-\pi/R^2})}{1-e^{-\pi/R^2}-\rho(\Delta,R)}.
$$
\end{corollary}
\textit{Proof:} First, observe that 
$$
\left\|P_{\Delta^c}(H-V_\varphi\sigma(H))\right\|_1\leq  \|n\|_1\leq \varepsilon.
$$
Hence,
$$
\left\|V_\varphi\big(f-\sigma(H)\big)\right\|_1
$$
$$
=\left\|P_{\Delta^c}V_\varphi\big(f-\sigma(H)\big)\right\|_1+\left\|P_\Delta V_\varphi\big(b-\sigma(g)\big)\right\|_1
$$
$$
\leq \left\|P_{\Delta^c}\big(V_\varphi f-H\big)\right\|_1+\left\|P_{\Delta^c}\big(H-V_\varphi\sigma(H)\big)\right\|_1$$
$$
+\delta(\Delta)\left\|V_\varphi\big(f-\sigma(H)\big)\right\|_1
$$
$$
\leq 2\varepsilon+\delta(\Delta)\left\|V_\varphi\big(f-\sigma(H)\big)\right\|_1,
$$
which concludes the proof using \eqref{delta-bound} and \eqref{planarsparselessone}.
\medskip

There are other approaches to the recovery of sparse time-frequency
representations which concentrate on the set-up of finite sparse
time-frequency representations \cite{SparsityTF,STF2}.

 Another
consequence of Theorem \ref{main} is the following refined $L_{1}$ uncertainty
principle for the STFT (see \cite[Proposition 3.3.1]{Charly} and \cite{BDJ,GM,RT} for other uncertainty principles for the STFT). 

\begin{corollary}\label{uncertainty}
Suppose that $f\in M^{1}$ satisfies $\Vert V_{\varphi }f\Vert _{1}=1$ and
that $\Delta \subset \mathbb{R}^{2}$ and $\varepsilon \geq 0$ are such that 
\begin{equation*}
1-\varepsilon\leq\int_{\Delta }|V_{\varphi }f(x,\omega )|dxd\omega,
\end{equation*}%
then 
\begin{equation*}
1-\varepsilon \leq \inf_{R>0}\left( \frac{\rho (\Delta ,R)}{1-e^{-\pi
/R^{2}}}\right) \leq |\Delta |.
\end{equation*}
\end{corollary}
\medskip

In particular, Corollary \ref{uncertainty} shows that the mass of the STFT of a function cannot be concentrated on sets that are locally small over the whole time-frequency plane. 

Our arguments to prove Theorem \ref{main} are an adaptation of Selberg's argument for the large sieve (see 
\cite{AB,Bombieri,Montgomery}), along the lines of \cite{DonohoLogan}. The analysis
reveals that, at least in the continuous case, dealing with joint
time-frequency representations leads to considerable simplifications, due to
the existence of local reproducing formulas \cite{Seip0}. This is not
surprising since, as observed earlier by Daubechies \cite{Daubechies} and
Seip \cite{Seip0}, the study of joint time-frequency restriction operators
with a gaussian window tends to be simplified. In particular, the functions best concentrated in a disc have a simple
explicit formula when written in the phase space. This is in contrast with the
classical time and band-limiting problem which has been studied in detail by Landau \cite{Land} (see also \cite{AP} for an alternative approach).

\section{Modulation and Fock spaces}

We will follow notations and definitions from \cite{Charly}. The \emph{%
Bargmann transform} on $\mathbb{C}$ is defined by 
\begin{equation}
\mathcal{B}f(z)=2^{1/4}\int_{\mathbb{R}}f(t)e^{2\pi t\cdot z-\pi
t^{2}-\pi z^{2}/2}dt\text{.}  \label{Barg}
\end{equation}%
Writing $z=x+i\omega $, a simple calculation shows that 
\begin{equation}
e^{-i\pi x\omega}V_{\varphi }f(x,-\omega )
=\mathcal{B}f(z)e^{-\pi |z|^{2}/2}\text{.%
}  \label{Bargmann}
\end{equation}%
Let $\mathcal{F}_{p}(\mathbb{C})$ be the space of entire functions equipped
with the norm 
\begin{equation*}
\Vert F\Vert _{\mathcal{L}_{p}}^{p}:=\int_{\mathbb{C}}|F(z)|^{p}e^{-\pi
p|z|^{2}/2}dz\text{,}
\end{equation*}%
if $1\leq p<\infty $ and 
\begin{equation*}
\Vert F\Vert _{\mathcal{L}_{\infty }}:=\sup_{z\in \mathbb{C}}|F(z)|e^{-\pi
|z|^{2}/2}\text{,}
\end{equation*}%
if $p=\infty $. The Bargmann transform is a unitary operator from $L^{2}(%
\mathbb{R})$ to $\mathcal{F}_{2}(\mathbb{C})$ and extends to a bijective
operator from $M^{p}$ to $\mathcal{F}_{p}(\mathbb{C})$, for $p\geq 1$ (see 
\cite{ST}, or \cite{AbrGr} for a proof that extends to polyanalytic Fock spaces). As in \cite{S}, we define the translation operator $\mathcal{T}_{w}
$ on $\mathcal{F}_{p}(\mathbb{C})$ as follows: 
\begin{equation*}
\mathcal{T}_{w}F(z):=e^{\pi \overline{w}\cdot z-\pi |w|^{2}/2}F(z-w).
\end{equation*}%
It acts isometrically on every $\mathcal{F}_{p}(\mathbb{C})$, $p\geq 1$. The
corresponding convolution is 
\begin{align}
F\ast G(z)&:=\int_{\mathbb{C}}F(w)G(z-w)e^{\pi (z\cdot \overline{w}%
-|w|^{2})}dw\nonumber \\ &=\int_{\mathbb{C}}F(w)\mathcal{T}_{w}G(z)e^{-\pi |w|^{2}/2}dw%
\text{.}
\end{align}

\section{Proof sketch of main results}

\subsection{\textbf{Concentration estimates}}

We say that a function $G\in \mathcal{L}_{1}(\mathbb{C})$ is concentrated on 
$\Omega \subset \mathbb{C}$ if $\Vert (I-P_{\Omega })G\Vert _{\mathcal{L}%
_{1}}=0$. Our main results will follow from the following statement  which corresponds to \eqref{DLl1}. We will give a full proof and more general results in \cite{AS}.

\begin{proposition}
\label{major-prop} Suppose that there exists $G\in \mathcal{L}_{\infty }(\mathbb{C})$ which is concentrated on $\Omega $, such that $F\mapsto
F\ast G$, $\mathcal{F}_{1}(\mathbb{C})\rightarrow \mathcal{F}_{1}(\mathbb{C})
$ is bounded and boundedly invertible. Then 
\begin{align}\label{eq-norm-bound}
\int_{\mathbb{C}}|F|d\mu 
\leq \Vert G\Vert _{\mathcal{L}%
_{\infty }}\cdot \Lambda (\mu ,\Omega )\cdot\nu (G)\cdot  \|F\|_{\mathcal{L}_1}\text{,}  \nonumber
\end{align}
where 
$
\nu (G):=\sup_{\Phi \in \mathcal{F}_{1}(\mathbb{C})}\left( \dfrac{\Vert \Phi
\Vert _{\mathcal{L}_{1}}}{\Vert \Phi\ast G \Vert _{\mathcal{L}_{1}}}\right) 
$
and 
\begin{equation*}
\Lambda (\mu ,\Omega ):=\sup_{w\in \mathbb{C}}\left(\int_{w+\Omega}e^{\pi |z|^{2}/2}d\mu (z)\right)\text{.}
\end{equation*}
\end{proposition}

\textit{Proof sketch:} For  $F\in \mathcal{F}_{1}(\mathbb{C})$, there exists $F^{\ast }\in \mathcal{F}_{1}(\mathbb{C})$ unique such that $F=F^{\ast }\ast
G$. 
Hence, replacing $F$ by $F^\ast\ast G$ and using, one after another, $\Vert (I-P_{\Omega
})G\Vert _{\mathcal{L}_{1}}=0$, Fubini's theorem, and H{\"{o}}lder's inequality ($p=1$) yields
$$
\int_{\mathbb{C}}|F(z)|d\mu (z)\leq \Vert G\Vert _{\mathcal{L}_{\infty }}\cdot \Lambda (\mu ,\Omega )\cdot 
\Vert F^{\ast }\Vert _{\mathcal{L}_{1}}.
$$
The observation that $\dfrac{\Vert F\|_{\mathcal{L}_{1}}}{\|F^\ast \ast G\|_{\mathcal{L}_{1}}}=1$ thus implies our statement.
\subsection{Proof of Theorem 1}

Define 
\begin{equation*}
d\mu(z):=\chi_\Delta(z)e^{-\pi |z|^{2}/2}dz\text{,}
\end{equation*}
with $\Delta \subset \mathbb{C}$ some subset of nonzero measure. Consequently, $\Vert F\Vert _{\mathcal{L}_{1,\mu }}=\Vert P_{\Delta }F\Vert _{\mathcal{L}%
_{1}}$ and 
setting  $\Omega =D_{1/R}$ yields
\begin{equation*}
\Lambda (\mu ,D_{1/R} )
=\sup_{z\in \mathbb{C}}|\Delta \cap (z+D_{1/R} )|=\rho(\Delta,R).
\end{equation*}%
Let $F$ be entire and $R>0$, then for any $z\in \mathbb{C}$,
the following local reproducing formula holds \cite{S}: 
\begin{equation}
F(z)=(1-e^{-\pi /R^{2}})^{-1}\cdot(F\ast\chi_{D_{1/R}})(z)\text{.}
\label{local-reproducing}
\end{equation}%
Now, let $R>0$, choosing $G=G_{R}:=\chi _{D_{1/R}}$
yields that convolution with $G_{R}$ gives a bounded and invertible operator on $%
\mathcal{F}_{1}(\mathbb{C})$. 
Then $\Vert G_{R}\Vert _{\mathcal{L}_{\infty }}=1$, 
\begin{equation*}
\nu (G_{R})=\sup_{\Phi \in \mathcal{F}_{1}}\left( \dfrac{\Vert
\Phi \Vert _{\mathcal{L}_{1}}}{\Vert\Phi\ast G_{R} \Vert _{\mathcal{L}_{1}}}%
\right) =\frac{1}{1-e^{-\pi /R^{2}}}\text{,}
\end{equation*}%
and Proposition \ref{major-prop} yields 
\begin{equation*}
\frac{\Vert P_{\Delta }F\Vert _{\mathcal{L}_{1}}}{\|F\|_{\mathcal{L}_1}}
 \leq\frac{\rho (\Delta ,R)}{1-e^{-\pi /R^{2}}}\text{.}
\end{equation*}%
This proves the result for $F\in \mathcal{F}_{1}(\mathbb{C})$. Since the
Bargmann transform {extends to a bijective operator from $M^{1}$ to $%
\mathcal{F}_{1}(\mathbb{C)}$, there exists }$f\in ${$M^{1}$ such that}%
\begin{equation*}
F(z)=\mathcal{B}f(z)=e^{-i\pi x\omega +\pi |z|^{2}/2}V_{\varphi
}f(x,-\omega )\text{.}
\end{equation*}%
This completes the proof.

\section*{Acknowledgement}
L.D.~Abreu and M. Speckbacker were supported by the Austrian Science
Foundation (FWF) START-project FLAME (\textquotedblleft Frames and Linear
Operators for Acoustical Modeling and Parameter Estimation\textquotedblright
, Y 551-N13)



\begin{thebibliography}{99}
\bibitem{AbrGr} L. D. Abreu, K. Gr\"{o}chenig. Banach Gabor
frames with Hermite functions: polyanalytic spaces from\ the Heisenberg
group. \emph{Appl. Anal.}, 91: 1981-1997, 2012.

\bibitem{AD} L. D. Abreu, M. D\"{o}rfler. An inverse problem for
localization operators. \emph{Inverse Problems}, 28:  115001, 2012.

\bibitem{AGR} L.~D. {A}breu, K.~{G}r{\"{o}}chenig, J.~L. {R}omero. On
accumulated spectrograms. \emph{Trans. Amer. Math. Soc.}, 368: 3629-3649,
2016.

\bibitem{AP} L. D. Abreu, J. M. Pereira. Measures of localization and quantitative {N}yquist densities. \emph{Appl. Comput. Harmon. Anal.}, 38 (3): 524-534, 2015

\bibitem{AS} L. D. Abreu, M. Speckbacher. \emph{In preparation.}

\bibitem{AB} C. Aubel, H. B\"{o}lcskei. Vandermonde matrices with nodes in
the unit disk and the large sieve. \emph{arXiv:1701.02538 preprint}, 2017.


\bibitem{BDJ} A. Bonami, B. Demange, P. Jaming. Hermite functions and
uncertainty principles for the Fourier and the windowed Fourier transforms. 
\emph{Rev. Math. Iberoamericana}, 19: 23-55, 2003.

\bibitem{Bombieri} E. Bombieri.  \emph{Le grand crible dans la th\'{e}orie
analytique des nombres}. Soci\'{e}t\'{e} Math\'{e}mathique de France, 1974.

\bibitem{CRT} E. J. Cand\'{e}s, J. Romberg, T. Tao. Robust uncertainty
principles: Exact signal reconstruction from highly incomplete frequency
information. \emph{IEEE Trans. Inf. Theor.}, 52: 489-509, 2006.

\bibitem{Donoho} D. L. Donoho. Compressed Sensing. \emph{IEEE Trans. Inform. Theor.}, 52(4): 1289-1306, 2006.

\bibitem{DonohoStark} D. L. Donoho, P. B. Stark. Uncertainty principles and
signal recovery. \emph{SIAM J. Appl. Math.}, 49: 906-931, 1989.

\bibitem{DonohoLogan} D. L. Donoho, B. F. Logan. Signal recovery and the large
sieve. \emph{SIAM J. Appl. Math.}, 52: 577-591, 1992.

\bibitem{Daubechies} I.~{D}aubechies. Time-frequency localization operators:
a geometric phase space approach. \emph{IEEE Trans. Inform. Theor.}, 34:
605-612, 1988.

\bibitem{FeiModulation} H.~G. Feichtinger. \newblock Modulation spaces on
locally compact abelian groups. \newblock In \emph{Proceedings of
``International Conference on Wavelets and Applications" 2002}, pages
99--140, Chennai, India, 2003. \newblock Updated version of a technical
report, University of Vienna, 1983.

\bibitem{Fei} H. G. Feichtinger. On a new Segal algebra. \emph{Monatsh. Math.%
} 92(4): 269--289, 1981.

\bibitem{FG} H. G. Feichtinger, K. Gr\"{o}chenig. Banach spaces related to
integrable group representations and their atomic decompositions, I. \emph{%
J. Funct. Anal.}, 86(2): 307-340, 1989.

\bibitem{DeMarie} F. DeMari, H. G. Feichtinger, K. Nowak. Uniform
eigenvalue estimates for time-frequency localization operators, \emph{J.
London Math. Soc.\ }, 65: 720-732, 2002.

\bibitem{FGHKR} H. F\"{u}hr, K. Gr\"{o}chenig, A. Haimi, A. Klotz, J. L.
Romero. Density of sampling and interpolation in reproducing kernel
Hilbert spaces.  \emph{arXiv:1607.07803v2 preprint}, 2016.

\bibitem{GM} K. Gr\"{o}chenig, E. Malinnikova. Phase space localization of
Riesz basis for $L^{2}(\mathbb{R}^{d})$. \emph{Rev. Mat. Iberoam}, 29: 
1003-1031, 2013.

\bibitem{JamPow} P. Jaming, A. Powell. Time-frequency concentration of
generating systems. \emph{Proc. Amer. Math. Soc}, 139:  3279-3990, 2011.

\bibitem{JamPowJFA} P. Jaming, A. Powell. Uncertainty principles for
orthonormal sequences. \emph{J. Funct. Anal.}, 243: 611-630, 2007.

\bibitem{Ly} Y. Lyubarskii. Frames in the Bargmann space of entire
functions. \emph{Entire and subharmonic functions}, 167-180, Adv. Soviet
Math., 11, Amer. Math. Soc., Providence, RI (1992).

\bibitem{Eugenia} E. Malinnikova. Orthonormal sequences in $L^{2}(R^{d})$
and time frequency localization. \emph{J. Fourier Anal. Appl. }, 16: 
983-1006, 2010.

\bibitem{Montgomery} H. L. Montgomery. The analytic principle of the large
sieve. \emph{Bull. Amer. Math. Soc.}, 84: 547-567, 1978.

\bibitem{FR} S. Foucart, H. Rauhut. \emph{A Mathematical Introduction to
Compressed Sensing}. Springer, New York, 2013.

\bibitem{Charly} K. Gr\"{o}chenig, \emph{\ Foundations of Time-Frequency
Analysis}. Birkh\"{a}user, Boston, 2001.

\bibitem{Land} H. J. Landau. Necessary Density Conditions for Sampling and
Interpolation of Certain Entire Functions. \emph{Acta Math.}, 117: 
37-52, 1967.

\bibitem{Logan} B. F. Logan. \emph{Properties of high-pass signals}. Ph.D.
thesis, Columbia University, New York, 1965.

\bibitem{SparsityTF} G. E. Pfander, H. Rauhut. Sparsity in time-frequency
representations. \emph{J. Four. Anal. Appl.}, 16: 233-260, 2010.

\bibitem{STF2} G. E. Pfander, H. Rauhut, J. Tanner. Identification of
matrices having a sparse representation. \emph{IEEE Trans. Signal Proc.},
56: 5376-5388, 2008.

\bibitem{RS95} J. Ramanathan, T. Steger. Incompleteness of sparse
coherente states. \emph{Appl. Comp. Harm. Anal.}, 2: 148-153, 1995.

\bibitem{RT} B. Ricaud, B. Torresani. A survey of uncertainty
principles and some signal processing applications. \emph{Adv. Comp. Math.}, 40(3): 629-650, 2014.

\bibitem{Seip0} K. Seip. Reproducing formulas and double orthogonality in
Bargmann and Bergman spaces. \emph{SIAM J. Math. Anal.}, 22: 856-876, 1991.

\bibitem{S} K. Seip, R. Wallst\'{e}n. Density Theorems for sampling and
interpolation in the Bargmann-Fock space I, \emph{J. Reine Angew. Math.},
429: 91-106, 1992.

\bibitem{ST} M. Signahl, J. Toft. Mapping properties for the Bargmann
transform on modulation spaces. \emph{J. Pseudo Diff. Oper. Appl.}, 3: 1-30,
2012.

\bibitem{Tao} T. Tao. An uncertainty principle for cyclic groups of prime
order. \emph{Math. Res. Lett.}, 12: 121-127, 2005.
\end{thebibliography}
\end{document}